# Some generalizations of Kannan's theorems via $\sigma_c$-function

Suprokash Hazra[1★] and Satish Shukla[2]

[1]Ramanujan School of Mathematical Sciences (RSMS),

Pondicherry University, Puducherry, India -605014.

[2]Department of Applied Mathematics, Shri Vaishnav Institute of Technology & Science,

Sanwer Road, Gram Baroli, Indore (M.P.) India- 453331.

**Abstract**

In this article we go on to discuss about various proper extensions of Kannan's two different fixed point theorems, introducing the new concept of $\sigma_c$ function ; which is independent of the three notions of simulation function, manageable functions and R-functions. These results are the analogous to some well known theorems, and extends several known results in this literature.

**Keywords:** Fixed point, coincidence point, Kannan's mapping, Simulation function, *R*-function, manageable function, $\sigma_c$-function.

**Mathematics Subject Classification** : 2000, 47H10, 54H25

## 1 Introduction:

The fixed point theory is one of the most useful and essential tools of nonlinear analysis. Banach (see, [1]) has given the most important and fundamental theorem of this branch by defining the concept of contraction operators. After that, so many theorems and generalizations of it has been made over the course of time. Recently, Khojasteh et al. (see, [5]) introduced the notion of simulation function and Du and Khojasteh (see, [6]) presented a very close but independent concept of manageable function, both of which give a new way to extend the Banach's fixed point result. However, López de Hierro and Shahzad (see, [7]) has given the concept of R-function (generalized concept of both simulation and manageable function) to obtain the extension of Banach's theorem for the R-contraction operator. On the other hand, Kannan (see, [2, 3]) found a particular type of operators which are not necessarily continuous, but satisfies the fixed point property on complete metric spaces. The class of operators found by Kannan and that of Banach are independent of each other (see, [2, 3]). Here, in this article, we prove several proper generalizations of Kannan's theorems, by finding fixed points and coincidence points for two set of operators, via the new concept of $\sigma_c$-functions. These new generalizations also extends several known theorems like Koparde-Waghmode theorem (See [12]) and Patel-Deheri's theorem (See [13]) ; by finding the analogous results of Malceski theorem (See [11]).

---

E-mail addresses : hazrasuprokash@gmail.com (Suprokash Hazra);

satishmathematics@yahoo.co.in (Satish Shukla).

★Corresponding author.

## 2  Preliminaries :

Consider, T, S : (X, d) → (X,d) as the two operators on metric space (X,d). T is said to have a fixed point c in X, if T(c) = c; and is said to be have coincidence point c in X, of the pair (T,S), if Tc = Sc. The space X is said to satisfy the coincidence property with respect to the pair (T,S) if there is atleast one point c for which Tc = Sc. The iterates of the self-mapping T are the mappings $\{T^n : X \to X\}_{n \in \mathbb{N} \cup \{0\}}$ defined by $T^0 = Id_X$, and $T^{n+1} = T \circ T^n$ for all n∈N. Given a point $x_0 \in X$, the Picard sequence of T based on $x_0$ is the sequence $\{x_n\}_{n \in \mathbb{N} \cup \{0\}}$ given by $x_{n+1} = Tx_n$ for all n ∈N. Clearly, $x_n = T^n x_0$ for all n∈N. T is said to be asymptotically regular at point x ∈ X if $\lim_{n \to \infty} d(T^n x, T^{n+1} x) = 0$. T is said to be sequentially convergent if, for each sequence $\{x_n\}$ the following holds true: if $\{Tx_n\}$ converges, then $\{x_n\}$ also converges. We say a sequence $\{x_n\}_{n \geq 0}$ to be S-bounded if $\{Sx_n\}_{n \geq 0}$ is bounded and S-Cauchy if $\{Sx_n\}_{n \geq 0}$ is cauchy sequence.

We now state the Kannan's two theorems, of which we find the generalizations.

**Theorem 2.1** : (Kannan, [3]) : *If T be an operator on a complete metric space X, satisfying the condition that, $d(Tx, Ty) \leq \alpha [d(Tx, x) + d(Ty, y)]$, for all x, y ∈X, where 0 < α < 1/2, then T has unique fixed point in X.*

**Theorem 2.2** : (Kannan; [2]) *Let X be a metric space with d as metric. Let T be a map of X into itself such that, (i) $d(Tx, Ty) \leq \alpha [d(Tx, x) + d(Ty, y)]$ ∀x,y ∈X, where 0 < α <1/2, (ii) T is continuous at a point c ∈X; and, (iii) There exists a point p ∈X such that the sequence of iterates $\{T^n(p)\}$ has a subsequence $\{T^{n_i}(x)\}$ converging to c, Then c is the unique fixed point of T.*

We state some definitions; starting with the concept of Simulation functions, which was initiated by Khojasteh et al. (see [5]), to show a new way to study, fixed point theory.

**Definition 2.3** :(Simulation function, [5]) : Let $\zeta: [0, \infty) \times [0, \infty) \to R$ be a mapping, then ζ is said to be a simulation function if it satisfies the following conditions,

(ζ1) $\zeta(0,0) = 0$;

(ζ2) $\zeta(t,s) < s - t$ for all t,s > 0;

(ζ3) if $\{t_n\}, \{s_n\}$ are sequences in $(0, \infty)$ such that, if $\lim_{n \to \infty} t_n = \lim_{n \to \infty} s_n > 0$; then, $\limsup_{n \to \infty} \zeta(t_n, s_n) < 0$.

**Example 2.4** :Let $\zeta: [0, \infty) \times [0, \infty) \to R$, such that, $\zeta(t, s) = \psi(s) - \phi(t)$ for all t,s ∈ $[0, \infty)$, where $\phi, \psi: [0, \infty) \to [0, \infty)$ are the two continuous functions such that $\psi(t) = \phi(t) = 0$ if and only if t = 0 and $\psi(t) < t \leq \phi(t)$ for all t > 0. More examples of simulation function can be found on [5].

**Definition 2.5** : (Manageable function, [6]) : A function $\eta: R \times R \to R$ is said to be manageable if the following conditions holds:

(η1) $\eta(t,s) < s - t$ for all s,t >0

($\eta$2) For any bounded sequence $\{t_n\} \subset (0,\infty)$ and any non increasing sequence $\{s_n\} \subset (0,\infty)$, we have that, $\limsup_{n \to \infty} [\{t_n + \eta(t_n, s_n)\}/ s_n] < 1$.

Several examples of manageable functions can be found on [6].

**Definition 2.6** : (*R*-function, [7]) : For a non empty $A \subseteq R$, a function $\rho : A \times A \to R$ is said to be an *R*-function if it satisfies the following two conditions.

($\rho$1) If $\{a_n\} \subset (0,\infty) \cap A$ is a sequence such that $\rho(a_{n+1}, a_n) > 0$, $\forall n \in N$, then $\{a_n\} \to 0$

($\rho$2) If $\{a_n\}, \{b_n\} \subset (0,\infty) \cap A$ are two sequences converging to the same limit $L \geq 0$; satisfying that, $L < a_n$, and $\rho(a_n, b_n) > 0$, $\forall n \in N$, then $L = 0$.

In some cases, the following additional property is also considered,

**($)** If $\{a_n\}, \{b_n\} \subset (0,\infty) \cap A$ are two sequences such that $\{b_n\} \to 0$ and $\rho(a_n, b_n) > 0$, $\forall n \in N$, then $\{a_n\} \to 0$.

Various examples and properties of R-function can be found on [7].

**Remark 2.7** (See [7]) : Every simulation function and manageable function is an *R*-function that also satisfies the property **($).**

**Definition 2.8** : (Geraghty function, [8]) A Geraghty function is a function $\phi : [0,\infty) \to [0,1)$ such that if $\{t_n\} \subset [0,\infty)$ and $\{\phi(t_n)\} \to 1$, then $\{t_n\} \to 0$.

**Definition 2.9** : (L-function, [9]) : A function $\phi : [0,\infty) \to [0,\infty)$ will be called an L-function if ;

(a) $\phi(0) = 0$ ;

(b) $\phi(t) > 0$ for all $t > 0$, and

(c) for all $\varepsilon > 0$, there exists $\delta > 0$ such that $\phi(t) \leq \varepsilon$ for all $t \in [\varepsilon, \varepsilon + \delta]$.

**3 Main Results :**

Before going into our main theorems, we shall introduce some definitions, as follows :

**Definition 3.1** : ($\sigma_c$-function) : Let $A \subseteq R$ be a non empty set and $c (\geq 1)$ a real constant. we introduce a function $\sigma_c : A \times A \to R$, is said to be an $\sigma_c$-function, if it satisfies the following two conditions :

($\sigma$1) If $\{a_n\} \subset (0,\infty) \cap A$ is a sequence such that $\sigma_c(a_n, a_{n-1} + a_n) > 0$, $\forall n \in N$, then $\{a_n\} \to 0$

($\sigma$2) If $\{a_n\}, \{b_n\} \subset (0,\infty) \cap A$ are two convergent sequences such that, $cL = \lim b_n \geq \lim a_n = L \geq 0$; satisfying that, $\sigma_c(a_n, b_n) > 0$, $\forall n \in N$, then $L = 0$.

We denote by $\Sigma_c^A$ as the set of all $\sigma_c$ functions on A, and we write simply $\Sigma_c$ for $\Sigma_c^{[0,\infty)}$.

**Remark 3.2 :** Unlike Remark 2.7, every simulation function (or manageable function) is not a $\sigma_c$-function. Infact the notion of $\sigma_c$-function is completely independent, from the three notions mentioned in definitions 2.3, 2.5, and 2.6 can be seen by the following examples.

**Example 3.3 :** Let $\gamma: [0,\infty) \times [0,\infty) \to R$ be the function defined, for all $t,s \in [0,\infty)$, and for $c > 0$ any reals, by,

$\gamma(t,s) = s/2 - 3t/2$, when $t < s$,

$\qquad = 0$, when $t \geq s$.

We, show $\gamma$ is a $\sigma_c$-function on $[0, \infty)$ (for all $0 < c < 3$), which satisfy the condition **($)**.

($\sigma 1$) If $\{a_n\} \subset (0,\infty)$ is a sequence such that $\gamma(a_n, a_{n-1}+a_n) > 0$, then, $\forall\, n \in N$

$0 < \gamma(a_n, a_{n-1}+a_n) = (a_{n-1}+a_n)/2 - 3a_n/2 = a_{n-1}/2 - a_n \Rightarrow a_n < a_{n-1}/2$ and hence $\{a_n\} \to 0$.

($\sigma 2$) If $\{a_n\}, \{b_n\} \subset (0,\infty)$ are two convergent sequences such that, $cL = \lim b_n \geq \lim a_n = L \geq 0$, satisfying that, $\gamma(a_n, b_n) > 0, \forall\, n \in N$, then, $0 < \gamma(a_n, b_n) = b_n/2 - 3a_n/2 \Rightarrow 0 < 3a_n/2 < b_n/2 \Rightarrow 0 \leq 3L \leq cL < 3L \Rightarrow L=0$.

**($)** If $\{a_n\}, \{b_n\} \subset (0,\infty)$ are two sequences such that $\{b_n\} \to 0$ and $\gamma(a_n,b_n) > 0, \forall\, n \in N$, then,

$0 < \gamma(a_n, b_n) = b_n/2 - 3a_n/2 \Rightarrow 0 < 3a_n/2 < b_n/2 \Rightarrow \{a_n\} \to 0$ as $\{b_n\} \to 0$ is given.

But $\gamma$ not a simulation function (or manageable function), because by considering $a_n = b_n = 1\, \forall\, n \in N$, we get that, $\limsup\, \gamma(a_n,b_n) = 0$, not satisfying the condition ($\zeta 3$) (or ($\eta 1$)).

**Example 3.4 :** Let, $\beta : [0,\infty) \times [0,\infty) \to R$ be the function defined, for all $t, s \in [0,\infty)$, by, $\beta(t,s) = s/2 - t$. Then clearly $\beta$ is a simulation function. Infact it is both manageable and R-function. But $\beta$ does not satisfy the property ($\sigma 1$). Because if we take $a_n = 1+1/n$, and $b_n = a_n + a_{n-1}, \forall\, n \in N$, then $\beta(a_n, b_n) = (a_{n-1} - a_n)/2 > 0$ but $a_n$ does not converge to 0.

**Example 3.5 :** Consider $g : [0,\infty) \times [0,\infty) \to R$ be the function defined, for all $t,s \in [0,\infty)$, by,

$g(t,s) = -1$, if $t \leq s$

$\qquad = 1$, if $t > s$,

Then $\forall\, n \in N$, and for every $\{a_n\} \subset (0,\infty)$, $g(a_n, a_{n-1}+a_n) < 0$. Hence condition ($\sigma 1$) is vacuously true.

Also, if $\{a_n\}, \{b_n\} \subset (0,\infty)$ are two convergent sequences such that, $cL = \lim b_n \geq \lim a_n = L \geq 0$, (for $c > 1$) satisfying that, $g(a_n, b_n) > 0, \forall\, n \in N$, then, it implies that, $a_n > b_n\, \forall\, n \in N \Rightarrow L \geq cL$ (for $c > 1$) $\Rightarrow L = 0$. So g satisfies ($\sigma 2$) and hence, g is a $\sigma_c$-function (for $c > 1$).

But g is not $R$-function can be seen by considering $a_n = n, \forall\, n \in N$.

**Remark 3.6 :** The domain of the functions $\gamma$ and g can be chosen any subset A of R, rather than $[0,\infty)$, to get examples of $\sigma_c$-functions with different domain. Also the above function $\beta$ satisfies ($\zeta 2$) but still not a $\sigma_c$-function, which would be very important when we state our conditions for the extension of Kannan mapping to get fixed point.

We will now consider some more examples of functions using Geraghty functions and L-functions.

**Example 3.7 :** Consider a function, $\pi : [0,\infty) \to R$, such that, $\pi(t) \leq t$, $\forall\, t \in [0,\infty)$. Now define a new function as, $\Theta_\pi : [0,\infty) \times [0,\infty) \to R$ defined by, $\Theta_\pi(t,s) = \alpha\,\pi(s) - t$, for all $t, s \in [0,\infty)$, with, $0 < \alpha < 1/2$ is a $\sigma_c$-function for all $1 \leq c \leq 2$, which also satisfies the property **($)**, can be readily seen.

**Note 3.8 :** In, [7], it is highlighted about an important property of L-function $l$, that, $l(t) \leq t$, $\forall\, t \in [0,\infty)$. So by previous example, for every L-function $l$, a function defined by; $\Theta_l(t,s) = \alpha l(s) - t$, for all $t, s \in [0,\infty)$; $(0 < \alpha < 1/2)$ is a $\sigma_c$-function for all $1 \leq c \leq 2$, which also satisfies the property ($)

**Proposition 3.9 :** *Let, $g : [0,\infty) \to [0,1)$ be a Geraghty function. Define $\Theta_g : [0,\infty) \times [0,\infty) \to R$ defined by, $\Theta_g(t,s) = \alpha g(s)s - t$, for all $t, s \in [0,\infty)$, with, $0 < \alpha \leq 1/2$, is a $\sigma_c$-function for all $1 \leq c \leq 2$, which also satisfies the property ($).*

*Proof* : For $0 < \alpha < 1/2$, the proof is clear from Example 3.7 and the fact that, $g(s) < 1$ ie. $g(s)s < s$.

Now for $\alpha = 1/2$, we have, $\Theta_g(t,s) = g(s)s/2 - t$, for all $t, s \in [0, \infty)$.

($\sigma$1)  If $\{a_n\} \subset (0,\infty)$ is a sequence such that $\Theta_g(a_n, a_{n-1}+a_n) > 0$, then, $\forall\, n \in N$, we have,

$0 < g(a_{n-1}+a_n)(a_{n-1}+a_n)/2 - a_n \Rightarrow 0 < a_n < g(a_{n-1}+a_n)(a_{n-1}+a_n)/2 < (a_{n-1}+a_n)/2 \Rightarrow a_n < (a_{n-1}+a_n)/2 \Rightarrow a_n < a_{n-1}$

So, $\{a_n\}$ is strictly monotone decreasing sequence of positive reals hence convergent to L (say).

Hence, $L \leq \lim g(a_{n-1}+a_n).2L/2 \leq 2L/2 = L \Rightarrow \lim g(a_{n-1}+a_n) = 1$, and so by the property of Geraghty function we have that, $(a_{n-1}+a_n) \to 0$ as $n \to \infty$; which implies $a_n \to 0$ as $n \to \infty$.

($\sigma$2) If $\{a_n\}, \{b_n\} \subset (0,\infty)$ are two convergent sequences such that, $cL = \lim b_n \geq \lim a_n = L \geq 0$, satisfying that, $\Theta_g(a_n, a_{n-1}+a_n) > 0$, then, we have $g(b_n)b_n/2 - a_n > 0 \Rightarrow a_n \leq g(b_n)b_n/2 < b_n/2 \Rightarrow L \leq \lim g(b_n)cL/2 < cL/2$

$\Rightarrow \lim g(b_n) = 1$ and so, $b_n \to 0$ as $n \to \infty$, and so $a_n \to 0$ as $n \to \infty$.

By similar arguments we can check the property **($)** and this completes the proof of the proposition.

**Definition 3.10 :** Let $T : X \to X$ be an operator. A sequence $\{x_n\}_{n \geq 0}$ in X, is said to satisfy the asymptotic regularity property, with respect to T, if, $\lim_{n \to \infty} d(Tx_{n+1}, Tx_n) = 0$. Now, if $\{x_n\}_{n \geq 0}$ is formed by the picard's interation, ie. $x_n = T^n x_0$, then T is simply turns out to be the asymptotically regular at the base point $x_0$.

**Definition 3.11** (See [10]) : Given two self-mappings $T, S : X \to X$ and a sequence $\{x_n\}_{n \geq 0} \subseteq X$, we say that $\{x_n\}_{n \geq 0}$ is a Picard sequence of the pair (T,S) (based on $x_0$) if $Sx_{n+1} = Tx_n$ for all $n \geq 0$. We say that X verifies the CLR(T,S)-property, if there exists on X a Picard sequence of (T,S) based on some point $x_0$.

**Definition 3.12 :** Let $T : X \to X$ be an operator. A sequence is said to be S-asymptotically similar with respect to T, if $\lim_{n \to \infty} d(TSx_n, SSx_n) = 0$,

**Definition 3.13 :** Let (X,d) be a metric space, and S : X → X be a function. A mapping T: X → X is called a $\Sigma_c$-S-Kannan with respect to some $\sigma_c \in \Sigma_c$, if, it satisfy the condition that,

$$\sigma_c (d(Tx,Ty), d(Tx,Sx) + d(Ty,Sy)) > 0 \ \forall \ x,y \in X. \qquad \ldots (3.13.1)$$

And is called, $\Sigma_c$-Kannan with respect to some $\sigma_c \in \Sigma_c$, if, it satisfy the condition that,

$$\sigma_c (d(Tx,Ty), d(Tx,x) + d(Ty,y)) > 0 \ \forall \ x,y \in X. \qquad \ldots (3.13.2)$$

**Lemma 3.14 :** *Let (X,d) be a metric space verifies the CLR(T,S) -property and let T: X → X be a $\Sigma_c$-S-Kannan, with respect to some $\sigma_c \in \Sigma_c$, then Picard sequence of (T,S) (based on $x_0$) satisfies, either the coincidence property with respect to the pair (T,S); or, the asymptotically regularity property with respect to the operator T.*

*Proof* : Given that, (X,d) verifies the CLR(T,S) -property. So there exists on X a Picard sequence $\{x_n\}_{n \geq 0}$ of (T,S) based on some point $x_0$ of X. which satisfying the condition that, $Sx_{n+1} = Tx_n$ for all n ≥ 0.

Now we have the following two cases.

**Case I:** We assume that, $Tx_p = Tx_{p-1}$, for some p∈N.

Then $Tx_p = Tx_{p-1} = Sx_p$ and X satisfies the coincidence property with respect to the pair (T,S).

**Case II :** Now we assume, $Tx_n \neq Tx_{n-1}, \forall \ n \in N$.

Now as T:X→X be a $\Sigma_c$-S-Kannan, with respect to some $\sigma_c \in \Sigma_c$, we have that,

$$\sigma_c (d(Tx,Ty), d(Tx,Sx) + d(Ty,Sy)) > 0 \ \forall \ x,y \in X..$$

So we have, $\sigma_c (d(Tx_{n+1}, Tx_n), d(Tx_{n+1}, Sx_{n+1}) + d(Tx_n, Sx_n)) > 0, \forall \ n \in N$.

Which implies, $\sigma_c (d(Tx_{n+1}, Tx_n), d(Tx_{n+1}, Tx_n) + d(Tx_n, Tx_{n-1})) > 0$, [As, $Sx_{n+1} = Tx_n$ for all n ≥ 0] ; $\forall \ n \in N$.

Choose $a_n = d(Tx_{n+1}, Tx_n)$, then, $a_{n-1} = d(Tx_n, Tx_{n-1})$; and so, $a_n > 0$ with $\sigma_c (a_n, a_{n-1} + a_n) > 0, \ \forall \ n \in N$.

Then by ($\sigma$1) we get, $\{a_n\} \to 0$. This clearly shows that, Picard sequence of (T,S) (based on $x_0$) satisfies the asymptotically regularity property with respect to the operator T.

**Lemma 3.15 :** *Let (X,d) be a metric space verifies the CLR(T,S) -property ;and T: X → X be a $\Sigma_c$-S-Kannan with respect to some $\sigma_c \in \Sigma_c$, satisfying, either the condition that, $\sigma_c (t, s) < s – t$ for all s, t > 0, or satisfying the property* ($). *Then the Picard sequence $\{x_n\}$ of the pair (T,S) (based on $x_0$), is a S-bounded sequence.*

*Proof* : Assume that $\{x_n\}$ is not S-bounded and we prove by contradiction. Without loss of generality we assume that $Sx_{n+p} \neq Sx_n \ \forall \ n, p \in N$ ;and so clearly, $Tx_{n+p} \neq Tx_n \ \forall \ n, p \in N$. Now as $\{x_n\}$ is not S-bounded, for each k, there exists two subsequences $\{Sx_{n_k}\}$ and $\{Sx_{m_k}\}$ of $\{Sx_n\}$ with $k \leq n_k < m_k$ ,that, for each $k \in N$, $m_k$, $n_k$ are the minimum integers, such that,

$$d(Sx_{n_k}, Sx_{m_k}) > 1 \text{ and } d(Sx_{n_k}, Sx_p) \leq 1 \text{ for } n_k \leq p \leq m_k - 1. \quad \ldots (3.15.1)$$

**Case I:** Now, if T is a $\Sigma_c$-S-Kannan with respect to some $\sigma_c \in \Sigma_c$, satisfying the condition that,

$$\sigma_c(t, s) < s - t, \text{ for all } s, t > 0.$$

So we have that,

$d(Tx,Ty) < d(Tx,Sx) + d(Ty,Sy)$, for those $x,y \in X$, of which both the sides provides non-zero entry.

Now, since we have assumed, $Tx_{n_k-1} = Sx_{n_k} \neq Sx_{n_k-1}$, and $Tx_{m_k-1} = Sx_{m_k} \neq Sx_{m_k-1}$,

We have that, $d(Tx_{n_k-1}, Tx_{m_k-1}) < d(Tx_{n_k-1}, Sx_{n_k-1}) + d(Tx_{m_k-1}, Sx_{m_k-1})$. [As, $Tx_{n+p} \neq Tx_n$]

Now clearly, $1 < d(Sx_{n_k}, Sx_{m_k}) = d(Tx_{n_k-1}, Tx_{m_k-1}) < d(Tx_{n_k-1}, Sx_{n_k-1}) + d(Tx_{m_k-1}, Sx_{m_k-1})$. [As,T is $\Sigma_c$-S-Kannan]

This implies, $1 < d(Tx_{n_k-1}, Sx_{n_k-1}) + d(Tx_{m_k-1}, Sx_{m_k-1}) = d(Tx_{n_k-1}, Tx_{n_k-2}) + d(Tx_{m_k-1}, Tx_{m_k-2})$

Now we see that both the entries on the right hand side are the subsequence of $d(Tx_n, Sx_n) = d(Tx_n, Tx_{n-1})$, such that, $Tx_n \neq Tx_{n-1}$ for all $n \in N$.

Then by Case II of Lemma 3.14, then Picard sequence of (T,S) (based on $x_0$) satisfies the asymptotically regularity property with respect to the operator T.

So, taking limit on both sides as $k \to \infty$, we get, $1 \leq 0$, which is a contradiction.

Hence, $\{x_n\}$ is S-bounded.

**Case II :** Now suppose, $\sigma_c$ if satisfies the property **($).

Then, we choose, $a_k = d(Sx_{n_k}, Sx_{m_k})$ and $b_k = d(Tx_{n_k-1}, Sx_{n_k-1}) + d(Tx_{m_k-1}, Sx_{m_k-1})$.

Then by, the given condition, we have, $a_k, b_k > 0$ satisfying $\sigma_c(a_k, b_k) > 0$ with $b_k \to 0$.

So by the property **($)** we have $a_k \to 0$, which is again a contradiction to (3.15.1),

This proves the lemma.

Now we use a similar type of idea as given in, [5] to prove the next lemma,

**Lemma 3.16 :** *Let (X,d) be a metric space verifies the CLR(T,S) -property; and T: X $\to$ X be a $\Sigma_c$-S-Kannan with respect to some $\sigma_c \in \Sigma_c$, satisfying, <u>either</u> the condition that, $\sigma_c(t, s) < s - t$ for all $s,t > 0$, <u>or</u>, satisfying the property* ($) . *Then the Picard sequence $\{x_n\}$ of the pair (T,S) (based on $x_0$), is a S-cauchy sequence.*

*Proof* : Consider $C_n = \sup\{ d(Sx_i, Sx_j) : i,j \geq n\}$. Note that the sequence $\{C_n\}$ is a monotonically decreasing sequence of positive reals and by Lemma 3.15, the sequence $\{x_n\}$ is S-bounded, therefore $C_n < \infty \; \forall \; n \in N$.

Thus $\{C_n\}$ is monotone, bounded sequence, hence convergent.

So there exists $C \geq 0$ such that $\lim_{n \to \infty} C_n = C$.

Now, if $C > 0$ then by the definition of $C_n$, for every $k \in N$ there exists $n_k, m_k$ such that $m_k > n_k \geq k$ and

$$C - 1/k < d(Sx_{m_k}, Sx_{n_k}) \leq C,$$

Hence, $\lim_{k \to \infty} d(Sx_{m_k}, Sx_{n_k}) = C$.

**Case I :** Now, suppose T is a $\Sigma_c$-S-Kannan with respect to some $\sigma_c \in \Sigma_c$, satisfying the condition that,

$\sigma_c(t, s) < s - t$, for all s, t > 0.

So we have that,

$d(Sx_{n_k}, Sx_{m_k}) = d(Tx_{n_k-1}, Tx_{m_k-1}) < d(Tx_{n_k-1}, Sx_{n_k-1}) + d(Tx_{m_k-1}, Sx_{m_k-1})$

Which implies, $d(Sx_{n_k}, Sx_{m_k}) < d(Tx_{n_k-1}, Tx_{n_k-2}) + d(Tx_{m_k-1}, Tx_{m_k-2})$

So, by previous argument, taking limit on both sides as $k \to \infty$, we get,

$\lim_{k \to \infty} d(Sx_{m_k}, Sx_{n_k}) = C \leq 0$,

This is a contradiction to the fact that, C > 0. Hence C = 0.

**Case II :** Now if $\sigma_c$ satisfies the property **($)**.

Then, we choose, $a_k = d(Sx_{n_k}, Sx_{m_k})$ and $b_k = d(Tx_{n_k-1}, Sx_{n_k-1}) + d(Tx_{m_k-1}, Sx_{m_k-1})$.

Then by, the given condition, we have, $a_k, b_k > 0$ satisfying $\sigma_c(a_k, b_k) > 0$ with $b_k \to 0$

So by the property **($)** we have $a_k \to 0$.

Hence C = 0, and this completes the lemma.

Now we state our theorem for $\sigma_c$-function, which analogous to the Theorem 2.2 (Kannan; [2]).

**Theorem 3.17 :** *Let (X,d) be a metric space verifies the CLR(T,S) -property and let $\{x_n\}_{n \geq 0}$ is a Picard sequence of the pair (T,S) (based on $x_0$). Let T be a map of X into itself such that,*

*(i) T: X $\to$ X be a $\Sigma_c$-S-Kannan with respect to some $\sigma_c \in \Sigma_c$ (for c=1), satisfying <u>either</u> the condition that,*
*(a) $\sigma_c(t, s) < s - t$ for all s, t > 0, <u>or</u> (b) $\sigma_c$ satisfying the property ($).*
*(ii) T and S both continuous at a point $Sq \in X$ ; and,*
*(iii) The picard sequence $Tx_n$ has a subsequence $\{Tx_{n_k}\}$ converging to Sq ; Then the following is true ;*

<u>Either</u> *(A) There is a coincidence point of the pair (T,S); <u>or both</u>*
*(B) $\{x_n\}_{n \geq 0}$ is S-asymptotically similar, implies, Sq is a fixed point of T; <u>and</u>*
*(B*) Sq is a fixed point of S, implies, Sq is a fixed point of T.*

*Proof* : (A) If atleast two consecutive terms of the picards sequence of pair (T,S) are equal, then by the case I of the Lemma 3.14, T has a coincidence point, which proves (A).

(B) Now assume that no terms of the picards sequence of the pair (T,S) are equal.

Given that, T is continuous at Sq $\in$ X, $\{T(Tx_{n_k})\}$ converging to TSq. We assume, TSq $\neq$ Sq, and will arrive at a contradiction. As, TSq $\neq$ Sq, we consider two disjoint open balls, say $B(TSq, r_1)$ and $B(Sq, r_2)$, with centres at TSq, Sq, and radius $r_1, r_2$ respectively.

We choose, r = min $\{r_1, r_2, d(TSq, Sq)/4\} > 0$.

Now, as the subsequence $\{Tx_{n_k}\}$ converging to Sq, and $\{T(Tx_{n_k})\}$ converges to TSq; there exists a positive integer M, such that, for all k > M, we have that,

$$Tx_{n_k} \in B(Sq, r) \text{ and } T(Tx_{n_k}) \in B(TSq, r),$$

And so clearly, for each k > M, we have that,

$0 < 4r < d(TSq, Sq) \leq d(TSq, T(Tx_{n_k})) + d(Tx_{n_k}, T(Tx_{n_k})) + d(Tx_{n_k}, Sq) < 2r + d(Tx_{n_k}, T(Tx_{n_k}))$

That is, $d(Tx_{n_k}, T(Tx_{n_k})) > 2r > 0$ ..... (3.17.1)

**Case I**: Now, if the condition (*i*)-(a) is satisfied, then we have that,

$0 < d(Tx_{n_k}, T(Tx_{n_k})) < d(Tx_{n_k}, Sx_{n_k}) + d(T(Tx_{n_k}), S(Tx_{n_k})) = d(Tx_{n_k}, Tx_{n_k-1}) + d(TSx_{n_k+1}, SSx_{n_k+1})$

Now, by the given condition (B), $\{x_n\}_{n \geq 0}$ is S-asymptotically similar; and using Lemma 3.14, we see that, the right hand side tends to 0.

So we get, $\lim_{k \to \infty} d(Tx_{n_k}, T(Tx_{n_k})) = 0$, which is a contradiction to (3.17.1).

Hence, TSq = Sq, that is, Sq is a fixed point of T.

**Case II**: If the condition (*i*)-(*b*) is satisfied, then we have,

$\sigma_c(d(Tx_{n_k}, T(Tx_{n_k})), d(Tx_{n_k}, Sx_{n_k}) + d(T(Tx_{n_k}), S(Tx_{n_k}))) > 0$.

We assume, $a_k = d(Tx_{n_k}, T(Tx_{n_k}))$ and $b_k = d(Tx_{n_k}, Sx_{n_k}) + d(T(Tx_{n_k}), S(Tx_{n_k}))$.

Then, since $\lim_{k \to \infty} b_k = 0$. which implies $\lim_{k \to \infty} a_k = 0$

It is a contradiction to (3.17.1).

Hence, T Sq = Sq, i.e., Sq is a fixed point of T

(B*) Now assuming the condition (i)-(a) to be satisfied by, $\sigma_c$; i.e, $\sigma_c(t, s) < s - t$, for all s, t > 0,

So we have that,

$$\sigma_c(d(Tx_{n_k}, T(Tx_{n_k})), d(Tx_{n_k}, Sx_{n_k}) + d(T(Tx_{n_k}), S(Tx_{n_k}))) > 0,$$

Choose, $a_k = d(Tx_{n_k}, T(Tx_{n_k}))$ and $b_k = d(Tx_{n_k}, Sx_{n_k}) + d(T(Tx_{n_k}), S(Tx_{n_k}))$,

we see that $a_k \to d(Sq, TSq) = L$ and $b_k \to d(TSq, SSq) = d(TSq, Sq) = L$, [by the assumption made in (c)]

Now, since T: X → X be a $\sigma_c$-Kannan with respect to some $\sigma_c \in \Sigma_c$ (for c=1), and L= $\lim a_k \leq \lim b_k = L$ ; it then implies that L= 0 and hence, d(Sq, TSq) = 0.

So, TSq = Sq and Sq is a fixed point of T.

Assuming the condition (*i*)-(*b*) to be satisfied by, $\sigma_c$; i.e, if $\sigma_c$ satisfies the condition **($)**, then similarly one can obtain the results, by using Lemma 3.14 and Lemma 3.16. This completes the proof.

Now we state our theorem for $\sigma_c$-function, which analogous to the Theorem 2.1 (Kannan; [3]).

**Theorem 3.18 :** *Let (X, d) be a metric space, and verifies the CLR(T,S) -property and let T be a map of X into itself such that,*

*(i) (S(X),d) is a complete, (or (T(X),d) is complete),*

*(ii) T is a $\Sigma_c$-S-Kannan with respect to some $\sigma_c \in \Sigma_c$ (for c=1), satisfying, <u>either</u> of the condition that,*

*(a) $\sigma_c(t, s) < s - t$ for all s, t > 0, <u>or</u>, (b) satisfying the property ($).*

*then X satisfies the coincidence property with respect to the pair (T,S).*

*Proof* : Suppose, T is a $\Sigma_c$-S-Kannan with respect to some $\sigma_c \in \Sigma_c$ (for c=1), satisfying the condition,

$$\sigma_c(t, s) < s - t, \text{ for all } s, t > 0;$$

Hence, by Lemma 3.16, the Picard sequence $\{x_n\}$ of the pair (T,S) (based on $x_0$), is a S-cauchy sequence.

Now, since (X,d) is complete, $\{Sx_n\}$ is convergent and converges to a point say z (say) in X.

Then as, $Sx_{n+1} = Tx_n$ for all $n \geq 0$. $\{Tx_n\}$ is also convergent and converges to the same point z.

Also as $z \in S(X)$, there is atleast one point w (say) in X such that Sw = z.

Also, if T satisfy the condition (*ii*)-(a); i.e $\sigma_c(t, s) < t - s$ for all s,t > 0;

Then, $\sigma_c(d(Tx_n, Tw), d(Tx_n, Sx_n) + d(Tw, Sw)) > 0$,

Now choose $a_n = d(Tx_n, Tw)$ and $b_n = d(Tx_n, Sx_n) + d(Tw, Sw)$, $\forall n \in N$.

Then, $\lim_{n\to\infty} a_n = \lim_{n\to\infty} d(Tx_n, Tw) = \lim_{n\to\infty} d(Sx_{n+1}, Tw) = d(z, Tw) = d(Sw, Tw)$ ; and,

$\lim_{n\to\infty} b_n = \lim_{n\to\infty} [d(Tx_n, Sx_n) + d(Tw, Sw)] = \lim_{n\to\infty} [d(Sx_{n+1}, Sx_n) + d(Tw, Sw)] = d(Tw, Sw)$

Then clearly, $d(Tw, Sw) = \lim b_n \geq \lim a_n = d(Tw, Sw) = L \geq 0$, so L = 0. (by $\sigma$2)

Hence X satisfies the coincidence property with respect to the pair (T,S).

Assuming the $\sigma_c$- function satisfies the condition (*ii*)-(b) i.e the property **($)**, one can obtain the result similarly, using Lemma 3.14 and Lemma 3.16. This completes the proof.

**Corollary 3.19 :** *Let (X, d) be a complete metric space, and, let T be a map of X into itself such that, T verifies the condition that, $\sigma_c(d(Tx, Ty), d(Tx, x) + d(Ty, y)) > 0 \ \forall \ x, y \in X$; with respect to some $\sigma_c \in \Sigma_c$ (c=1), and satisfy <u>either</u>, the condition ($) <u>or</u>, that, $\sigma_c(t, s) < s - t$, for all s,t > 0; then T has* unique *fixed point in X.*

*Proof :* Putting, S = $Id_X$, i.e, S(x) = x for all $x \in X$ in theorem 3.18 we obtain the existence of the fixed point of T. Only to prove the uniqueness. Note that, for S = $Id_X$ the Picard sequence of the pair (T,S) based on some point $x_0 \in X$, now reduces into the Picard sequence $\{x_n\}$ of T based on $x_0$ for arbitrary $x_0 \in X$, where $x_n = T^n x_0$; and $\{x_n\}$ converges to u such that u is a fixed point of T. If possible, assume Tv = v with u ≠ v, for some $v \in X$. Then, we can choose a subsequence $\{x_{n_k}\}$ of $\{x_n\}$ such that $x_{n_k} \neq x_{n_k+1}$ and, $Tx_{n_k} \neq Tv$ for all $k \in N$ (because, if not so, then taking limits u = Tv = v).

Now if $\sigma_c(t, s) < s - t$, for all s, t > 0; we have that, $d(Tx_{n_k},Tv) < d(Tx_{n_k}, x_{n_k}) + d(Tv, v)$, and taking limits we get, $d(u, Tv) = 0$. i.e. $u = Tv = v$.

Also if, the condition **($)** is satisfied, then we choose,

$a_k = d(Tx_{n_k}, Tv)$ and $b_k = d(Tx_{n_k}, x_{n_k}) + d(Tv, v)$; and as $b_k \to 0$ as $k \to \infty$, we get $a_k \to 0$ as $k \to \infty$.

So, $d(u, Tv) = 0$. i.e. $u = Tv = v$. This completes the proof of the corollary.

**Corollary 3.20 :** $T: R^n \to R^n$ be a function such that $T(\underline{a}) = -T(-\underline{a})$, for some $\underline{a} \in R^n$, and satisfies the condition that, $||Tx - Ty|| \le 2\alpha (||Tx|| + ||Ty||)$, for a fixed $0 < \alpha < 1/2$. Then there is root of the function in $R^n$, i.e there is a point $\underline{c}$ in $R^n$ for which $f(\underline{c}) = \underline{0}$

*Proof :* Put $S(x) = -T(x)$ in the theorem 3.18. and consider $\sigma_c(t, s) = \alpha s - t$. $\forall$ $t, s \in [0, \infty)$, $0 < \alpha < 1/2$. (For the proof of being $\sigma_c$-function, see Remark 3.21, below)

**Remark 3.21 :** The theorems 3.17 and 3.18 are the generalizations of the theorem 2.2 (Kannan; [2]) and theorem 2.1 (Kannan; [3]) respectively. Because for a fixed $0 < \alpha < 1/2$ if we choose a function $\chi^\alpha$, as,

$\chi^\alpha : [0, \infty) \times [0, \infty) \to R$ defined by, $\chi^\alpha(t,s) = \alpha s - t$. $\forall$ $t, s \in [0, \infty)$,

then $\chi^\alpha$ is $\sigma_c$-function (for $0 < c < 2$), for every, $0 < \alpha < 1/2$.

We will show this fact.

($\sigma 1$) If $\{a_n\} \subset (0, \infty)$ is a sequence such that $\chi^\alpha(a_n, b_n) > 0$, where, $b_n = a_{n-1} + a_n$, then, $\forall$ $n \in N$

$0 < \chi^\alpha(a_n, b_n) = \alpha b_n - a_n = \alpha(a_{n-1} + a_n) - a_n = \alpha a_{n-1} - (1 - \alpha) a_n \Rightarrow a_n < \alpha a_{n-1}/(1 - \alpha)$,

and as $0 < \alpha < 1/2$ the quantity $\alpha/(1 - \alpha) < 1$ and hence $\{a_n\} \to 0$.

($\sigma 2$) If $\{a_n\}, \{b_n\} \subset (0, \infty)$ are two convergent sequences such that, $cL = \lim b_n \ge \lim a_n = L \ge 0$, satisfying that, $\chi^\alpha(a_n, b_n) > 0$, $\forall$ $n \in N$, then, $0 < \chi^\alpha(a_n, b_n) = \alpha b_n - a_n \Rightarrow 0 < a_n < \alpha b_n \Rightarrow 0 \le L \le c\alpha L < cL/2 < L \Rightarrow L = 0$. (as $c < 2$)

This proves that $\chi^\alpha$ is $\sigma_c$-function (for $c < 2$), for every, $0 < \alpha < 1/2$.

Now, if an operator satisfies the condition (3.13.1) of the definition 3.13; then for $\sigma_c = \chi^\alpha$ and $S = Id_X$ we have that,

$$\chi^\alpha(d(Tx,Ty), d(Tx, x) + d(Ty, y)) > 0, \ x, y \in X$$

$\Rightarrow d(Tx, Ty) \le \alpha[d(Tx,x) + d(Ty,y)]$, $\forall$ $x, y \in X$, ($0 < \alpha < 1/2$).

Which is the Kannan's contraction condition. All the other conditions of Theorem 3.17 and Theorem 3.18 can be reduced into desired form, easily.

**Remark 3.22 :** The identity function is the simplest operator having fixed point, without satisfying the Kannan's contraction (or not even Banach's contraction condition). But it does satisfy our $\Sigma_c$-S-Kannan contraction condition. Also our $\Sigma_c$-S-Kannan contraction condition is a proper extension of the Kannan's contraction condition even when we take S to be the identity operator. Also the theorem 3.17 and Theorem 3.18 are the two proper generalizations of Kannan's theorem. These all facts follows from the next few examples.

**Example 3.23 :** First let us consider the example 3.5, with a slight modification, that,

Consider $\omega : [0,\infty) \times [0,\infty) \to R$ be the function defined, for all $t,s \in [0,\infty)$, by,

$\omega(t,s) = -1$,  if $t < s$

$\phantom{\omega(t,s)} = 1$,  if $s \leq t$,

Then, in view of Example 3.5 we have $\omega \in \sum_c$ with c > 1. Now if we take our operator T to be the identity function then it is a $\sigma_c$-Kannan mapping. Since for any $x,y \in X$, we have that,

$\omega(d(Tx,Ty), d(Tx,x) + d(Ty,y)) = \omega(d(x,y), 0) = 1 > 0$.

But $\omega$ not a Kannan's mapping; because if so, then for $x \neq y$, we would have that,

$d(x,y) = d(Tx,Ty) \leq \alpha [d(Tx,x) + d(Ty,y)] = 0$, which is a contradiction.

Also if we take, X= any space, and T: $X \to X$ be any functions with atleast two different images, and we choose S to be any constant function, then X satisfy the CLR(T,S) property and T satisfying the $\Sigma_c$-S-Kannan contraction condition with respect to the above function $\omega$. Also if would have chosen the last considered T, pre-assuming that there exists a point c, that T(c)=c then T is satisfied all the conditions above without the condition (ii) of Theorem 3.18; and this shows that the condition (ii) of the Theorem 3.18, is sufficient but not necessary to have fixed point.

**Example 3.24 :** To prove our $\Sigma_c$-S-Kannan contraction condition is a proper extension of the Kannan's contraction condition, we consider X= {1,2,3,4,5} and define T: $X \to X$ as T(x) = 3 if $x \neq 4$ and T(4) = 2; and define S: $X \to X$ as, S(x) = 3, $x \neq 4$ and S(4) =5 . Then clearly T does not satisfy the Kannan's contraction can be seen by considering two points 3 and 4. But this T is satisfying the $\Sigma_c$-S-Kannan contraction condition with respect to the $\sigma_c \in \Sigma_c$, defined by $\sigma_c(t,s) = 3s/7 - t$ (for 0< c < 2). (By Remark 3.21). Infact X satisfying the CLR(T,S) property and S(X) is complete as well. So it satisfy all the property of the Theorem 3.18. and hence X satisfies the coincidence property with respect to the pair (T,S) and x=1,2,3 are the coincidence points of it.

**Remark 3.25 :** In [7] it was shown that every simulation function and every manageable function is an R-function that also satisfies the property **($)**. Notice that, the property ($\sigma$1) of $\sigma_c$ functions is quite different from the property ($\rho$1) of manageable functions. This property is a major difference between the manageable functions and the σc functions, and it plays an important role in the existence of fixed point of the class of σc-mappings. The following example shows that one cannot replace the property (σ1) by (ρ1) even when the property **($)** is satisfied.

**Example 3.26 :** Let X = {1,2,3} and consider the metric d(x,y) = |x−y| for all x,y ∈ X. Clearly (X,d) is complete. Now, consider a function τ : [0,∞) × [0,∞) → R be the function defined by τ(t, s) = 2s/3 − t, for all t,s ∈ [0,∞). Clearly τ is a manageable function which satisfies the property **($)**.

We now define an operator T : X → X by T(1) = T(2) = 3 and T(3) = 2.

To show τ(d(Tx,Ty),d(Tx,x) + d(Ty,y)) > 0 we show that d(Tx,Ty) ≤ 2[d(Tx,x) + d(Ty,y)]/3.

Then, 0 =d(T(1),T(2)) < 2 [d(T(1),1) + d(T(2),2)]/3;

1 =d(T(1),T(3)) < 2[d(T(1),1) + d(T(3),3)]/3 = 2(2+1)/3 =2,

1=d(T(2),T(3)) < 2[d(T(2),2) + d(T(3),3)]/3 = 2 [1 +1]/3 = 4/3.

Also, it is easy to see that τ satisfies the property (σ2) with 1 ≤ c < 3/2. Thus, all the conditions of Theorem 3.18, except the condition (σ1), are satisfied, and T has no fixed point.

We now move to find another different types of generalizations of Kannan's theorems, starting with the following definition, as follows.

**Definition 3.27** : In a metric space (X,d),  T: X→X and S: X→X be the two operators. Then T is said to be S-dominated $\Sigma_c$-Kannan mapping of degree w, if, with respect to some $\sigma_c \in$  $_c$,the following holds,

$\sigma_c$ (d$^w$(STx,STy), d$^w$(Sx,STx) + d$^w$(Sy,STy)) >0 , for all x,y ∈X, for any fixed w∈N.          ….. (3.27.1)

Now, In [11] Malceski proved the following generalization of Kannan' theorem:

**Theorem 3.28** : *Let (X,d) be a complete metric space, T: X→X and S: X→X  be a mapping such that it is continuous, injection and sequentially convergent. If α > 0,γ ≥ 0 and 2α+γ < 1, and , satisfies the condition that, d(STx,STy) ≤ α[d(Sx,STx)+ d(Sy,STy)] + γd(Sx,Sy); for all x,y∈X, then there is a unique fixed point of T.*

Now we find another different generalization of Theorem 2.1 (Kannan, [3]), by using $\sigma_c$ function,which is quite analogous to the above Theorem 3.28.

**Theorem 3.29**: *Let (X,d) is a metric space, and  T: X→X and S: X→X be the two operators; such that, T is S-dominated $\Sigma_c$-Kannan mapping of degree w, with respect to some $\sigma_c \in \Sigma_c$ (c=1), with the following conditions holds;*
*(1) S(X) is complete,*
*(2) S is injective;*
*(3) <u>Either</u> $\sigma_c$ (t,s) < s - t ;<u>or,</u> $\sigma_c$  satisfy the condition ($),*
*Then there exists a unique fixed point of the operator T, in X.*

*Proof:* Let {$x_n$}$_{n≥0}$ is formed by the picard's interation, ie. $x_n$=T$^n$$x_0$, with initial base point $x_0$. we assume, $x_{n+1}$ ≠ $x_n$ for any n; because if not so for some index k, then T($x_k$)= $x_{k+1}$ = $x_k$ and we would get fixed  point



and our proof would be over. So assume $x_{n+1} \neq x_n$ for all n. Now, by the injectivity of S implies that, for such a sequence, $Sx_{n+1} \neq Sx_n$ for all n.

Given, $\sigma_c (d^w(STx,STy), d^w(Sx,STx)+ d^w(Sy,STy)) > 0 \ \forall \ x,y \in X$.

So we have, $\sigma_c (d^w(STx_{n+1}, STx_n), d^w(STx_{n+1}, Sx_{n+1}) + d^w(STx_n, Sx_n)) > 0, \forall n \in N$.

Which implies, $\sigma_c (d^w(Sx_{n+2}, Sx_{n+1}), d^w(Sx_{n+2}, Sx_{n+1}) + d^w(Sx_{n+1}, Sx_n)) > 0$, [As, $x_{n+1} = Tx_n$ for all $n \geq 0$]

Choose $a_n = d^w(Sx_{n+2}, Sx_{n+1}) (>0)$, then, $a_{n-1} = d^w(Sx_{n+1}, Sx_n)$; and so $\sigma_c (a_n, b_n) > 0$, where, $b_n = a_{n-1} + a_n, \forall n \in N$.

So by ($\sigma$1) we have $d^w(Sx_{n+1}, Sx_n) = a_{n-1} \to 0$, as $n \to \infty$. That is $d(Sx_{n+1}, Sx_n) \to 0$ as $n \to \infty$. ...... (3.30.1)

**Claim**: $\{x_n\}$ is S-bounded.

*Proof of Claim*: Using the similar previous arguments, we assume, that $\{x_n\}$ is not S-bounded and we prove by contradiction. Without loss of generality we assume that $Sx_{n+p} \neq Sx_n \ \forall \ n, p \in N$; and so clearly, $Tx_{n+p} \neq Tx_n \ \forall \ n, p \in N$. Now as $\{x_n\}$ is not S-bounded, for each k, there exists two subsequences $\{Sx_{n_k}\}$ and $\{Sx_{m_k}\}$ of $\{Sx_n\}$ with $k \leq n_k < m_k$, that, for each $k \in N$, $m_k$ is the minimum integer such that,

$$d(STx_{n_k-1}, STx_{m_k-1}) = d(Sx_{n_k}, Sx_{m_k}) > 1 \text{ and } d(Sx_{n_k}, Sx_p) \leq 1 \text{ for } n_k \leq p \leq m_k \quad ..... (3.30.2)$$

**Subcase I**: Now if $\sigma_c(t,s) < s - t$ then, by (3.30.2), we have that,

$1 < d^w(Sx_{n_k}, Sx_{m_k}) = d^w(STx_{n_k-1}, STx_{m_k-1}) < d^w(STx_{n_k-1}, Sx_{n_k-1}) + d^w(STx_{m_k-1}, Sx_{m_k-1})$

$$= d^w(Sx_{n_k}, Sx_{n_k-1}) + d^w(Sx_{m_k}, Sx_{m_k-1})$$

Then by (3.30.1), taking limit on both sides as $k \to \infty$, we arrive at a contradiction, as previous.

**Subcase II**: Also if $\sigma_c$ satisfy the condition **($)**, then we assume,

$a_k = d^w(Sx_{n_k}, Sx_{m_k})$ and $b_k = d^w(Sx_{n_k}, Sx_{n_k-1}) + d^w(Sx_{m_k}, Sx_{m_k-1})$,

Then by, the given condition, we have, $\sigma_c(a_k, b_k) > 0$ with $b_k \to 0$ so by the property (S) we have $a_k \to 0$.

**Claim**: $\{x_n\}$ is S-Cauchy.

*Proof of Claim*: Likewise previous case, we consider $C_n = \sup\{ d(Sx_i, Sx_j) : i, j \geq n\}$.

Then by the S-boundedness and monotonicity of $C_n$ implies, $\lim_{n \to \infty} C_n = C$. (for some C). And as previous, if $C > 0$, then there exists $q_k, p_k$ with $p_k > q_k \geq k$, that $\lim_{k \to \infty} d(Sx_{p_k}, Sx_{q_k}) = C$. Then like previous cases (or by Subcase I, subcase II, for $q_k, p_k$) we proved the claim.

Now by condition (1) there exists a point z (say) in S(X), for which $Sx_n \to z$. Also as $z \in S(X)$, there is a point say $c \in X$ that $S(c) = z$.

Now as previous we consider, $a_n = d^w(STx_n, STc)$ and $b_n = d^w(STx_n, Sx_n) + d^w(STc, Sc)$, for all $n \in N$, such that, $a_n \to d^w(z, STc)$ and $b_n \to d^w(Sc, STc) = d^w(Sc, STc)$. So by ($\sigma$2) we have, $d^w(z, STc) = 0$, ie. $STc = z = Sc$ and, the injectivity of S shows that, $Tc = c$. Uniqueness of fixed point follows from the previous arguments.



**Remark 3.30** : Here we have only assumed the injective condition of S, and haven't use the continuity and sequential convergence criteria, but still the theorem remains true, if we assume that image of X under S is complete.

**Corollary 3.31** :(**Koparde-Waghmode theorem ;** See [12]): *If (X,d) is a complete metric space and, T:X→X is mapping such that, there exists, α∈(0,1/2) so that, the following is satisfied ,*

$$d^2(Tx,Ty) \leq \alpha [\ d^2(x,Tx) + d^2(y,Ty)], \text{ for all } x, y \in X,$$

*Then T would have unique fixed point in X.*

Proof : We assume our $\sigma_c$ function as $\sigma_c(t,s)$ = αs - t, with α∈(0,1/2), and consider our operator as, S(x)=x, and put w= 2, in the Theorem 3.29.

**Corollary 3.32 :** (**Patel-Deheri theorem** ; See [13]): *In a complete metric space, if the two operators S,T satisfy, d(STx, STy) ≤ α [ d(Sx, STx) + d(Sy, STy)], where 0 < α < 1/2 and S is continuous, injection and sequentially convergent; then T has a unique fixed point in X.*

*Proof* : We assume our $\sigma_c$ function as $\sigma_c(t,s)$ = αs - t, with α∈(0,1/2), and put w= 1, in the Theorem 3.29.

Now we prove another theorem, analogous to Theorem 2.2 (Kannan [2]) ,

**Theorem 3.33** : *Let (X,d) be a metric space. Let T,S be a map of X into itself such that,* T is S-dominated $\Sigma_c$ -Kannan mapping of degree w=1, with respect to some $\sigma_c$ function (c=1), with the following conditions holds;

(i) Either $\sigma_c(t, s) < s - t$ for all s,t > 0, <u>or</u>, $\sigma_c$ satisfying the condition ($)
(ii) There exists a point p such that, the picard sequence $Tx_n$ has a subsequence $\{Tx_{n_k}\}$ converging to q;
(iii) T and S both continuous at the point Sq∈X; and,
(iv) S is injective and continuous at Tq;
Then, q is the unique fixed point of T.

*Proof*: Suppose Tq≠q. Then, by (iv), STq≠Sq. We consider two disjoint open balls, say $B$(TSq, $r_1$) and $B$(Sq, $r_2$), with centres at Tq, q, and radius $r_1$, $r_2$ respectively and choose, 0< r < min {$r_1$, $r_2$, d(TSq, Sq)/3 }.

Now, as the subsequence $\{Tx_{n_k}\}$ converging to q, and S is continuous at both q and Tq. So $\{STx_{n_k}\}$ converges to Sq; and $\{STTx_{n_k}\}$ converges to STq. there exists a positive integer M, such that, for all k > M, we have that,

$STx_{n_k} \in B(Sq, r)$ and $STTx_{n_k} \in B(STq, r)$,

And so clearly, for each k > M , we have that,
Then, By (i), 0 < 3r < d(STq, Sq) ≤ d(STq, STTx_{n_k}) + d( STTx_{n_k}, STx_{n_k}) + d(STx_{n_k}, Sq) < 2r + d( STx_{n_k}, STTx_{n_k})



This implies that,

$$d(STx_{n_k}, STTx_{n_k}) > r. \quad\quad\quad ....(3.33.1)$$

Now, T is S-dominated $\Sigma_c$-Kannan mapping of degree w=1, so by using (*i*) and, using the similar ideas used in Theorem 3.17, S is asymptotically regular and hence we have contradiction to (3.33.1). And so STq= Sq or equivalently Tq=q. The uniqueness follows from the case is as previous.

**Example 3.34**: Also our S-dominated $\Sigma_c$-Kannan contraction condition is a proper extension of the Kannan's contraction condition we consider the example given in [10], [12] with a little brief.

X= {0, 1/4, 1/5,1/6,...} endowed with the Euclidean metric. Define T: X→X, by T(0)= 0, and T(1/n)= 1/(n+1). And consider, S :X→X to be S(0)=0 and S(1/n) = $1/n^n$, then T does not satisfy the Kannan's contraction condition for any constant >0, but it is S-dominated $\Sigma_c$-Kannan mapping of degree w=1 with respect to the $\sigma_c$-function defined by $\sigma_c(t,s)$ = s/3 – t (for 0 < c < 2). (By Remark 3.20). Also by defining S to be, S(0)=0 and S(1/n) = $1/[e^{2n}]$, T is a S-dominated $\Sigma_c$-Kannan with respect to the $\sigma_c$-function defined by,
$\sigma_c(t,s)$ = s/6– t (for 0 < c < 2). (By Remark 3.21).

We consider one more new example of this type, but in much simpler form.

**Example 3.35:** Consider the example discussed in Example 3.24. Then clearly T does not satisfy the Kannan's contraction can be seen by considering two points 3 and 4. But this T is satisfying the S-dominated $\Sigma_c$-Kannan, with respect to the $\sigma_c$-function defined by $\sigma_c(t,s)$ = $\alpha$s– t (for 0 < c < 2) ;for all fixed α∈(0,1/2).

## 4  Conclusion :

Here in this article we tried to find several extension of Kannan's two different fixed point results, by introducing the concept of $\Sigma_c$-S-Kannan operator and S-dominated $\Sigma_c$-S-Kannan operator of degree w; via the new concept of $\sigma_c$-function (shown to be independent of the other three concepts of Simulation function, Manageable function and R-function). These new generalizations also extends several known theorems in this branch and the similar ideas could be profitably extended to the three dimensional case of $\sigma_c$-functions and would be helpful to find the extension of the Fisher type of mapping and the similar type of operators. Now, by the aboves discussions, following interesting problems arise.

**Problem 1:** Can the Theorem 3.17 (or Theorem 3.18) and Theorem 3.29 (or Theorem 3.33) be proved with an operator satisfying the condition (3.13.1) (or condition 3.27.1) for a dense subset of X, instead of whole of X?

**Problem 2**: If we define an analogous definition 3.13, and Definition 3.27 as follows :

An operator T on a metric spaces is called $\Sigma_c$-S- Fisher if it satisfy, $\sigma_c$ (d(Tx,Ty), d(Tx,Sy) + d(Ty,Sx)) > 0 and said to be S-dominated $\Sigma_c$-Fisher mapping of degree w, if, with respect to some $\sigma_c$ function, if $\sigma_c$ (d$^w$(STx, STy), d$^w$(Sy,STx) + d$^w$(Sx,STy)) > 0 , for all x,y $\in$ X, for any fixed w $\in$ N. Then, can the Theorem 3.17 (or Theorem 3.18) and Theorem 3.29 (or Theorem 3.33) be proved, with an operator satisfying the above proposed conditions for the Fisher type of operator?

**Author's Contribution :**

The 1st author initiated and originated the whole work and contributed mainly and significantly to this article. The 2nd author also contributed with meticulous experiments, scintillating suggestions and corrections, after making the final manuscript of the paper. Both the authors read and approved the final manuscript of the paper.

**Conflicts of Interest :**

Authors declare that there is no conflicts of interest regarding the publication of this paper.

**Acknowledgements:**

The 1st author like to thank his teachers Prof. Dr. Biswajit Mitra, Dr. Anandamoy Mukhopadhyay and Dr. Kusal Chattopadhyay for being the constant source of inspiration and encouragement throughout the every part of his academic life. The same author also like to thank the 2nd author for taking a kind initiative in some context and valuable response to this work.